\documentclass [12pt,oneside] {amsart}
\usepackage{latexsym,amsfonts,amssymb,mathrsfs}
\usepackage [all]{xy}
\SelectTips{cm}{12}
\title{Representation and character theory in $2$-categories}
\author{Nora Ganter and Mikhail Kapranov}
\thanks{Ganter was supported by NSF grant DMS-0504539. Kapranov was
  supported by NSF grant DMS-0500565. The paper was completed while
  Ganter was visiting MSRI}
\date{\today}

\theoremstyle{plain}
\newtheorem {Thm} {Theorem}[section]
\newtheorem* {Thm*} {Theorem}
\newtheorem* {Prop*} {Proposition}
\newtheorem {Lem}[Thm] {Lemma}
\newtheorem {Prop}[Thm] {Proposition}
\newtheorem {Cor}[Thm] {Corollary}

\newtheorem {Obs}[Thm] {Observation}

\theoremstyle {definition}
\newtheorem {Def}[Thm] {Definition}

\newtheorem {Not}[Thm] {Notation}
\newtheorem {Rem}[Thm] {Remark}
\newtheorem {Exa}[Thm] {Example}
\newenvironment{Pf}[1]{{\noindent\sc Proof #1:}}{\qed\\}
\newcommand {\specialmap} [4] {\text {$ #1\negmedspace : #2 #3 #4 $}}
\newcommand {\map} [3] {\specialmap {#1} {#2}{\to} {#3}}
\newcommand {\isomap} [3] {\specialmap {#1} {#2}{\overset {\cong}
            {\longrightarrow}} {#3}}
\newcommand {\tisomap} [3] {\specialmap {#1} {#2}{\overset {\cong}
            {\Longrightarrow}} {#3}}

\newcommand {\Aut}{\operatorname{Aut}}
\newcommand {\id} {\operatorname{id}}

\newcommand {\at}[1] {\arrowvert_{#1}}

\newcommand {\ind} {\operatorname{ind}}
\renewcommand {\(} {\left(}
\renewcommand {\)} {\right)}

\newcommand {\sub} {\subseteq}

\newcommand {\CC} {\mathbb C}
\def\Pfeil{{\unitlength 1em\begin{picture}(0,.1)
\put(0,.1){\vector(1,0){1.5}}\end{picture}}}
\def\colimname{{\unitlength.1em
\raisebox{-2.7\unitlength}{\begin{picture}(15.5,9.5)(0,0)
\put(0,2.7){$\operatorname{lim}$}
\put(.05,-.1){\Pfeil}
\end{picture}}}}
\def\Colim{\mathop{\colimname}}
\newcommand {\on}[1] {\operatorname{#1}}
\newcommand {\BG} {{\on{BG}}}

\newcommand {\mC} {{\mathcal C}}
\newcommand {\mD} {{\mathcal D}}
\newcommand {\Cl}[1] {{\on{Cl}(#1;k)}}

\newcommand {\con}[2] {{#1}^{-1}{#2}{#1}}

\newcommand {\GG}[1] {{\frac{1}{|#1|}}}

\renewcommand {\ind}[2] {\on{ind}\arrowvert_{#1}^{#2}}
\newcommand {\inv} {^{-1}}

\newcommand {\LG} {\Lambda(G)}

\newcommand {\mmod} {/\!\!/}
\newcommand {\mor} {{\on{Mor}}}

\newcommand {\mR} {{\mathcal R}}

\newcommand {\MU} {{\on{MU}}}

\newcommand {\mV} {{\mathcal V}}

\newcommand {\ob} {{\on{Ob}}}

\newcommand {\omor} {$1$-morphism }
\newcommand {\omors} {$1$-morphisms }

\newcommand {\orb} {{\on{orb}}}

\newcommand {\pt} {\on{pt}}

\newcommand {\Repk}[1] {{\mathcal{R}ep(#1)}}
\newcommand {\res}[2] {\on{res}\arrowvert^{#1}_{#2}}
\newcommand {\RG} {R(G)}
\newcommand {\Set} {\on{Set}}

\newcommand {\Sn} {{\Sigma_n}}

\newcommand {\tensor}{\otimes}

\newcommand {\tcat} {$2$-category }
\newcommand {\tcats} {$2$-categories }
\newcommand {\tchar} {categorical character }
\newcommand {\thom} {\on{2-Hom}}

\newcommand {\tmors} {$2$-morphisms }
\newcommand {\trep} {$2$-representation }
\newcommand {\treps} {$2$-representations }
\newcommand {\tvect} {{\ensuremath \on{2-Vect}_k}}
\newcommand {\ttr}[1] {{\mathbb{T}r(#1)}}

\newcommand {\tr} {\on{Trace}}
\newcommand {\Vect} {\on{Vect}}
\newcommand {\Xind} {{\chi_{\on{ind}}}}

\newcommand {\ZZ} {\mathbb Z}

\begin{document}

\begin{abstract}
We develop a (2-)categorical generalization of the theory of
group
representations and characters. We categorify the concept
of the trace of a linear transformation, associating
to  any endofunctor of any small  category a set called its categorical
trace. In a linear situation, the categorical trace is a vector space
and we associate to any two commuting self-equivalences a number called
their joint trace. For a group acting on a linear category $\mathcal {V}$
we define an analog of the character which is the function
on commuting pairs of group elements given by the joint traces
of the corresponding functors. We call this function the 2-character of 
$\mathcal {V}$. 
Such functions of commuting pairs (and more generally, $n$-tuples)
 appear in the work of Hopkins, Kuhn and
Ravenel \cite{Hopkins:Kuhn:Ravenel} on equivariant Morava E-theories. We define the
concept of induced categorical representation and show that
the corresponding 2-character is given by the same formula as was obtained in 
\cite{Hopkins:Kuhn:Ravenel} for the transfer map in the second equivariant
Morava E-theory.

\end{abstract}
\maketitle
\pagestyle{plain}
\section{Introduction}
The  goal of this paper is to develop a (2-)categorical generalization
of the theory of group representations and characters.
It is classical that a representation $\varrho$ of a group $G$ is
often determined by its character
$$\chi(g)=tr(\varrho(g)),$$
which is a class function on $G$.

\vskip .3cm

Remarkably,
generalizations of
character theory turn up naturally in the context of homotopy theory.
Since this so called {\em Hopkins-Kuhn-Ravenel character theory}
motivated much of our work, we will start with a short and very informal
summary of it and some other homotopy theoretic topics.
Fix a prime $p$,
let $n$ be a natural number, and let $BG$
denote the classifying space 
of $G$. Assume that $G$ is finite. 
In \cite{Hopkins:Kuhn:Ravenel}, Hopkins, Kuhn and Ravenel computed the
ring $E^*_n(BG)$, where $E_n$ is a generalized cohomology theory
(depending on $p$), which
was introduced by Morava \cite{Rezk}.
The first Morava $E$-theory is $p$-completed $K$-theory,
 $$E_1=K{}_p\widehat{}.$$
Hopkins, Kuhn and Ravenel found
that elements $\chi\in E^*_n(BG)$ are most naturally described as
$n$-class functions, i.e.
functions
$$\chi(g_1,\dots,g_n)$$
defined on $n$-tuples of commuting elements of $G$ and invariant under
simultaneous conjugation. 
In the context of \cite{Hopkins:Kuhn:Ravenel}, all the $g_i$ are
required to have $p$-power order.
Hopkins, Kuhn and Ravenel actually make a much stronger case that the
$n$-class functions occurring in this way should be viewed as
generalized group characters: 
if $\alpha\negmedspace :H\hookrightarrow G$ is the
inclusion of a subgroup, then there is a map
$$
  \map{B\alpha}{BH}{BG},
$$
and in the stable homotopy category, one has a transfer map $\tau\alpha$
in the other direction. These maps make
the correspondence $G\mapsto E_n(BG)$ into a Mackey functor.
The map $E_n^*(\tau\alpha)$ sends a generalized character of $H$ to a
generalized character of $G$. If we stick with the analogy to
classical character theory, it plays the role of the induced 
representation. Hopkins, Kuhn and Ravenel compute its effect on
generalized characters. They find that it is described by the formula
\begin{equation}
  \label{hkr-ind-Eqn}
  E_n^*(\tau\alpha)(\chi)(g_1,\dots ,g_n) = \frac1{|H|}\sum_{\stackrel{s\in
      G \medspace\mid}{s\inv\underline g s\in H^n}}  
    \chi(s\inv g_1s,\dots,s\inv g_n s),  
\end{equation}
where $\underline g=(g_1,\dots,g_n)$ is an $n$-tuple of commuting
elements in $G$.
For $n=1$, this is the formula for the character of the induced
representation, cf.  \cite{Serre:representations}.

\vskip .3cm


The number $n$ is often referred to as the {\em chromatic level} of
the theory, see \cite{Ravenel} for general background on the
chromatic picture in homotopy theory. In the case $n=2$, the 
theory $E_2$ is an example of an elliptic cohomology
theory. For background on elliptic cohomology, we refer the
reader to the introduction of \cite{AHS}).

Just as the representation ring $R(G)$ may be viewed as equivariant
$K$-theory of the one point space, the ring $E_2(BG)$ is 
Borel equivariant $E_2$-theory of the one point space.
Elliptic cohomology is a field at the
intersection of several areas of mathematics, and there is a
variety
of fields that have motivated definitions of equivariant elliptic
cohomology. 
To name a few, we have 
Devoto's definition, motivated by
orbifold string theory
\cite{Devoto}, we have Grojnowski's work \cite{Grojnowski}, motivated
by the theory 
of loop groups, we have the axiomatic approach in \cite{GKV}, 
involving principal bundles over elliptic curves, we have a connection
to generalized Moonshine (cf.\ \cite{Ganter}, \cite{Devoto}, and
\cite{Baker:Thomas}), 
and we have recent constructions of Lurie and Gepner
\cite{Lurie}, \cite{Gepner}, 
which
satisfy axioms similar to those of \cite{GKV}
but formulated in the context of derived algebraic geometry.
Lurie's construction naturally involves $2$-groups.
It is remarkable that each of these constructions, in one way or
another, leads to class-functions on pairs of commuting elements of
the group. 

\vskip .3cm

What is lacking in  the above approaches is an analog of the notion
of representation which would produce the generalized characters by
means of some kind of trace construction.
In this paper, we supply such a notion (for $n=2$). 
It turns out that the right object to consider is
an action of $G$ on a category instead of a vector space or, more
generally, on an object of a $2$-category.

Our main construction is the so-called {\em categorical trace} of
a functor $A: \mathcal{V}\to\mathcal{V}$  from a small category  $\mathcal{V}$ 
to itself (or, more generally, of a
1-endomorphism of an object of a 2-category). This categorical trace
is a set, denoted  $\mathbb{T}r(A)$, see Definition \ref{trace}.
 If $k$ is a field, and the category  $\mathcal{V}$ is $k$-linear, then 
$\mathbb{T}r(A)$ is a $k$-vector space. In the latter case,
 given
two commuting self-equivalences $A,B: \mathcal{V}\to\mathcal{V}$,
we define their {\em joint trace} $\tau(A,B)$ to be the ordinary trace of
the linear transformation induced by $B$ on the vector space $\mathbb{T}r(A)$.
Here we assume that $\mathbb{T}r(A)$ is finite-dimensional. 
 
For a group $G$ acting on $\mathcal{V}$ this gives a 2-class function
called the 2-character of the categorical representation $\mathcal{V}$. 
Among other things, we define an induction
procedure
for categorical representations and show that it produces
the map (\ref{hkr-ind-Eqn}) on the level of characters.

\vskip .3cm

This  very simple and natural construction
 ties in with other geometric approaches to elliptic cohomology.
Already  Segal, in his Bourbaki talk \cite{Segal}, proposed
to look  for some kind of ``elliptic objects'' which are related
to vector bundles in the same way as 2-categories are related to
ordinary categories. While vector bundles can be equipped with
connections and thus with the concept of parallel transport along
paths, one expects elliptic objects to allow parallel transport along
2-dimensional membranes. 
Similarly, more 
recent
works \cite{Teichner:Stolz}, \cite{Baas:Dundas:Rognes},
\cite{Hu:Kriz}
aiming at geometric definitions of elliptic cohomology,
all  involve $2$-categorical constructions.

\vskip .3cm

The present paper is kept at an  elementary level and
does not require any knowledge of homotopy theory
(except for the final  section \ref{comp-HKR}
devoted to the comparison with  
 \cite{Hopkins:Kuhn:Ravenel}). Nevertheless, the connection
with homotopy theory and, in particular, with equivariant
elliptic cohomology was important for us. It provided us with
a motivation as well as with 
 a possible future field of applications. 

We also do not attempt to discuss actions of groups on
$n$-categories for $n>1$ which seem to be the right way to get
$(n+1)$-class functions.

\vskip .3cm

We have recently learned of a work in progress by Bruce Bartlett and
Simon Willerton \cite{Bartlett:Willerton} 
where, interestingly, the concept of the categorical trace also
 appears although
the motivation is different.

Our inspiration for this project came
from conversations with Haynes Miller. 
We would like to thank
Jim Stasheff and
 Simon Willerton
 for their remarks on  earlier
versions of this paper. We are grateful to the
referee, whose comments greatly helped to 
make the paper accessible for a larger audience. 
We would also like to thank
Matthew Ando, Alex Ghitza, and
Charles Rezk  for many helpful conversations. 

\section{Background on $2$-categories}
\label{2background-Sec}
\subsection{} Recall \cite{Maclane}  that a 2-category $\mathcal {C}$
consists of the following data:
\label{2cat-Sec}
\begin{enumerate}

\item{} A class of objects $\mathcal {O}b \, \mathcal{C}$.

\item For any $x,y\in \mathcal {O}b\,\mathcal{ C}$ a category $Hom_{\mathcal {C}}(x,y)$.
Its objects are called 1-morphisms from $x$ to $y$
 (notation $A: x\rightarrow y$).
We will also use the notation $1Hom_{\mathcal {C}}(x,y)$ for
 $\mathcal {O}b \,  Hom_{\mathcal C}(x,y)$. For any
 $A,B\in 1Hom_{\mathcal C}(x,y)$
morphisms from $A$ to $B$ in $Hom_{\mathcal {C}}(x,y)$ are called 2-morphisms
from $A$ to $B$ (notation $\phi: A\Rightarrow B$). We denote the set of such
morphisms by

$$2Hom_\mC(A, B) = Hom_{Hom_{\mathcal {C}}(x,y)}(A, B).$$

The composition in the category $Hom_{\mathcal{C}}(x,y)$ will be denoted by
$\circ_1$ and called the vertical composition. Thus if $\phi: A\Rightarrow B$
and
$\psi: B\Rightarrow C$, then $\psi\circ_1 \phi: A\Rightarrow C$.

\item The composition bifunctor
\begin{eqnarray*}
  Hom_{\mathcal{C}}(y,z) \times Hom_{\mathcal{C}}(x,y) & \rightarrow &
  Hom_{\mathcal{C}}(x,z) \\
  (A, B)&\mapsto &A\circ_0 B.
\end{eqnarray*}
In particular, for a pair of $1$-morphisms $A, B: x\rightarrow y$, a $2$-morphisms $\phi:
A\Rightarrow B$ between them, and a pair of $1$-morphisms $C, D:
y\rightarrow z$ 
and a $2$-morphism $\psi: C\Rightarrow D$ between them, there is a
composition
$\psi\circ_0 \phi: C\circ_0 A\Rightarrow D\circ_0 B$.

\item The natural associativity isomorphism
$$\alpha_{A,B,C}: (A\circ_0 B)\circ_0 C\Rightarrow A\circ_0(B\circ_0
  C).$$
It is given for any three composable 1-morphisms $A,B,C$ and satisfies the
  pentagonal
axiom, see \cite{Maclane}.

\item For any $x\in \mathcal{O}b\, \mathcal{C}$, a 1-morphism ${1}_x\in 1Hom_{\mathcal {C}}
(x,x)$ called the unit morphism, which comes equipped with $2$-isomorphisms

$$\epsilon_\phi: {1}_x\circ_0 A \Rightarrow A, \quad\text{for any } A: y\rightarrow x,$$
$$\zeta_\psi: B\circ {1}_x\Rightarrow B, \quad\text{for any } B: x\rightarrow z,$$

\hskip -1.7cm  satisfying the axioms of \cite{Maclane}.

\end{enumerate}

We denote by $1\mor\mC$ and $2\mor\mC$ the classes of all 1- and 2-morphisms of $\mC$.
We say that $\mathcal C$ is strict if all the  $\alpha_{A,B,C}$,
$\epsilon_\phi$ and $ \zeta_\psi$ are identities (in particular the
source
and target of each of them are equal). It is a theorem of Mac Lane and
Par\'e   that every 2-category can be replaced by a (2-equivalent)
strict one. See \cite{Maclane} for details.

\subsection{Examples}\label{2cat-Exa-Sec}
(a) The 2-category $\mathcal{C}at$ has, as objects, all small
categories, as morphisms their functors and as 2-morphisms natural transformations
of functors. This 2-category is strict.
 We will use the notation $Fun(\mathcal{A, B})$ for the set of functors
between categories $\mathcal A$ and $\mathcal B$ (i.e., 1-morphisms in $\mathcal{C}at$)
and $NT(F,\Phi)$ for the set of natural transformations between functors $F$ and $\Phi$.
Most of the examples of 2-categories can be embedded into $\mathcal {C}at$:
a 2-category $\mathcal {C}$ can be realized as consisting of categories with
some extra structure.

\vskip .3cm

(b) Let $k$ be a field. The 2-category $2Vect_k$, see \cite {Kapranov:Voevodsky}
 has, as objects, symbols
$[n]$, $n=0,1,2, ...$ For any two such objects $[m], [n]$ the category
$Hom([m], [n])$ has, as objects, 2-matrices of size $m$ by $n$, i.e., matrices
of vector spaces $A= \|A_{ij}\|$, $i=1, ..., m$, $j=1, ..., n$. Morphisms between
2-matrices $A$ and $B$ of the same size are collections of linear maps $\phi =
\{ \phi_{ij}: A_{ij}\rightarrow B_{ij}\}$. Composition of 2-matrices is given by the
formula
$$(A\circ B)_{ij} = \bigoplus_l A_{il} \otimes B_{lj}.$$
This 2-category is not strict. An explicit strict replacement was
constructed in \cite{Elgueta}.

\vskip .3cm

(c) The 2-category $\mathcal{B}im$ has, as objects, associative rings. If $R, S$
are two such rings, then $Hom_{\mathcal{B}im}(R, S)$ is the category of $(R,S)$-bimodules.
The composition bifunctor
$$Hom(S,T)\times Hom(R,S) \rightarrow Hom(R,T)$$
is given by the tensor product:
$$(M, N)\mapsto N\otimes_S M.$$
This 2-category is also not strict.

{ Relation to $\mathcal {C}at$:} To a ring $R$ we associate the category
$\on{Mod-R}$
of right $R$-modules. Then each $(R,S)$-bimodule $M$ defines a functor
$$\on{Mod -}R \rightarrow \on{Mod -}S, \quad P\mapsto P\otimes_R M.$$

The 2-category $2Vect_k$ is realized inside $\mathcal{B}im$ by associating to
$[m]$ the ring $k^{\oplus m}$. An $m$ by $n$ 2-matrix is the same as a
$(k^{\oplus m}, k^{\oplus n})$-bimodule.

We will denote by $\mathcal {B}im_k$ the sub-2-category in $\mathcal{B}im$
formed by $k$-algebras as objects and the same 1- and 2-morphisms as in
$\mathcal{B}im$.

\vskip .3cm

(d) Let $X$ be a CW-complex. The Poincare 2-category $\Pi(X)$ has, as objects,
points of $X$, as 1-morphisms Moore paths $[0,t]\rightarrow X$ and as
2-morphisms homotopy classes of homotopies between Moore paths.

\vskip .3cm

We will occasionally use the concept of a (strong) $2$-functor
$\map\Phi\mC\mD$ between \tcats $\mC$ and $\mD$. Such a $2$-functor
consists of maps $\ob\mC\to\ob\mD$, $1\mor\mC\to 1\mor\mD$ and
$2\mor\mC\to2\mor\mD$ preserving the composition of \tmors and
preserving the composition of \omors up to natural
$2$-isomorphisms. See \cite{Maclane} for details.

\subsection {2-categories with extra structure.}
We recall the definition of an enriched category from \cite{Maclane}
 or \cite{Kelly}.
Let $(\mathcal{A},\tensor,S)$ be a closed symmetric monoidal category,
so $\otimes$ is the monoidal operation and $S$ is a unit object.
\begin{Def}
  A category enriched over $\mathcal{A}$ (or simply a $\mathcal{A}$-category) $C$
 is defined in the same way as a category, with
  the morphism sets replaced by objects $\hom(X,Y)$ of $\mathcal{A}$ and
  composition replaced by $\mathcal{A}$-morphisms
  $$
    \hom(X,Y)\tensor\hom(Y,Z)\to\hom(X,Z),
  $$
  with units
  $$
    \map{1_X}S{\hom(X,X)},
  $$
  such that the usual associativity and unit diagrams commute.
\end{Def}

\begin{Exa}
Categories enriched over the category of abelian groups are
commonly
known as pre-additive categories. By an additive category one
means a pre-additive category possessing finite direct sums.
If $k$ is a field, categories enriched over the category of $k$-vector
spaces are known as $k$-linear categories.
\end{Exa}
\begin{Exa}
A strict \tcat is  the same as  a
category enriched over the category of small categories
with $\otimes$ being the direct product of categories (cf. \cite{Street}).
\end{Exa}
\begin{Def}\label{enriched-Def}
  Let $\mathcal{A}$ be a category. A {\em strict \tcat $\mathcal{C}$ enriched over
 $\mathcal{A}$}, or
  shorter an {\em $\mathcal{A}$-$2$-category}, is a category
  enriched over the category of small $\mathcal{A}$-categories.
\end{Def}
\begin{Def}
We define a strict pre-additive \tcat to be a \tcat enriched over the
  category of abelian groups. Let $k$ be a field. Then
a (strict)  $k$-linear \tcat is defined to be a
  \tcat enriched over the category of $k$-vector spaces.
Weak additive and $k$-linear 2-categories are defined in a similar way.
\end{Def}

We will freely use the concept of a triangulated category
\cite{Neeman},
\cite{Gelfand:Manin}. If $\mathcal {D}$
is triangulated, then we denote by $X[i]$ the $i$-fold iterated shift (suspension)
of an object $X$ in $\mathcal{D}$. We will denote
  $$Hom^\bullet_{\mathcal {D}}(X, Y) = \bigoplus_i Hom_{\mathcal {D}}(X, Y[i]).$$

We will call a 2-category $\mathcal {C}$ triangular if each
$Hom_{\mathcal{C}}(x,y)$ is made into a triangulated category and the composition
functor is exact in each variable.

\subsection {Examples.}\label{Examples} (a) The 2-category $\mathcal{B}im$
 is additive.
The 2-categories $2Vect_k$ and $\mathcal {B}im_k$ are $k$-linear.

\vskip .3cm

(b) Define the 2-category $\mathcal{DB}im$ to have the same objects as $\mathcal{B}im$,
i.e., associative rings. The category $Hom_{\mathcal{DB}im}(R,S)$
is defined to be the  derived category of complexes of $(R,S)$-bimodules
bounded above. The composition is given by the derived tensor product:

$$(M, N)\mapsto N\otimes_S^L M.$$

This gives a triangular 2-category.

\vskip .3cm

(c) The 2-category $\mathcal{V}ar_k$ has as objects smooth projective
algebraic varieties over $k$. If $X,Y$ are two such varieties, then
$$Hom_{\mathcal{V}ar_k}(X,Y) = D^b \mathcal{C}oh(X\times Y)$$
is the bounded derived category of coherent sheaves on $X\times Y$.
If $\mathcal {K}\in D^b\mathcal{C}oh(Y\times Z)$ and $\mathcal{L}\in
D^b\mathcal{C}oh(X\times Y)$, then their composition is defined by the
derived convolution

$$\mathcal{K} * \mathcal {L} = Rp_{13*} (p_{12}^*\mathcal {L}
\otimes^L p_{23}^*\mathcal{K}),$$
where $p_{12}, p_{13}, p_{23}$ are the projections of $X\times Y\times Z$ to the
products of two factors.  This again gives a triangular 2-category.

\vskip .2cm

Relation to $\mathcal{C}at$: To every variety $X$ we associate the category
$D^b\mathcal{C}oh(X)$. Then every sheaf $\mathcal {K}\in D^b\mathcal{C}oh(X\times Y)$
(``kernel'')
defines a functor
$$F_{\mathcal {K}}: D^b\mathcal{C}oh(X)\rightarrow D^b\mathcal{C}oh(Y), \quad
\mathcal F\mapsto Rp_{2*} (p_1^* \mathcal{F}\otimes^L \mathcal {K}),$$
and $F_{\mathcal{K}*\mathcal{L}}$ is naturally isomorphic to $F_{\mathcal{K}}
\circ F_{\mathcal{L}}$. It is not known, however,  whether the natural map
$$Hom_{D^b\mathcal{C}oh (X,Y)}(\mathcal{K}, \mathcal{L})
\to NT (F_{\mathcal{K}}, F_{\mathcal{L}})$$
is a bijection for arbitrary $\mathcal{K},\mathcal{L}$. So in practice
the source of this map is used as a substitute for its target.

\vskip .3cm

(d) The 2-category ${\Bbb R}{\mathcal A}n_k$ has, as objects, real analytic manifolds.
For any two such manifolds $X,Y$ the category $Hom_{\mathcal{CW}}(X,Y)$
is defined to be $D^b \mathcal{C}onstr (X\times Y)$, the bounded derived
category of ($\mathbb {R}$-) constructible sheaves of $k$-vector spaces on $X\times Y$,
see \cite{Kashiwara:Schapira}, Sect. 8.4.,   for background on constructible sheaves. 
The composition is defined
similarly to the above, with $p_{ij}^*$ understood as sheaf-theoretic direct
images rather than  $\mathcal{O}$-module-theoretic direct images.
This is a triangular 2-category.

\vskip .2cm

Relation to $\mathcal{C}at$: To every real analytic manifold $X$ we associate the 
category $D^b \mathcal{C}onstr (X)$. Then, as in (c), any sheaf 
$$\mathcal {K}\in D^b\mathcal{C}onstr(X\times Y)$$ can be considered
as a ``kernel'' 
defining a functor $$D^b \mathcal{C}onstr (X) \to D^b \mathcal{C}onstr (Y).$$

\vskip .3cm

(e) Let $\mathcal{A}b$ denote the category of all abelian categories. For any
such categories $\mathcal{A,B}$ the category $Fun(\mathcal{A,B})$ is again
abelian: a sequence of functors is exact if it takes any object into an exact
sequence. So we have a triangular 2-category $\mathcal{DA}b$ with same objects
as $\mathcal{A}b$ but $Hom(\mathcal{A,B}) = D^b  Fun(\mathcal{A,B})$.

\vskip .2cm

\section{The categorical trace}
\subsection{The main definition}
As motivation consider the $2$-category $2Vect_k$. In this situation
there is a  na\"\i ve
way to define the
``trace'' of a $1$-automorphism, namely as direct sum of the diagonal
entries of the matrix. This na\"\i ve notion of trace is equivalent
to the following definition that makes sense in any $2$-category $\mC$:
\begin{Def}\label{trace}
  Let $\mathcal{C}$ be a 2-category, $x$ an object of $\mathcal{C}$
  and $A: x\rightarrow x$ a 1-endomorphism of $x$. The categorical trace of $A$
  is defined as
  $$\ttr{A} = 2Hom_{\mathcal C} ({1}_x, A).$$
  If $\mathcal{C}$ is triangular,
  we write
  $$\mathbb{T}r^i(A) = \mathbb{T}r(A[i]), i\in\mathbb{Z}, \quad \mathbb{T}r^\bullet(A) =
  \bigoplus_{i} \mathbb{T}r^i(A).$$
\end{Def}
\begin{Rem}[Functoriality] Note that for each $x$, the categorical
  trace defines a functor
  \begin{eqnarray*}
    \mathbb{T}r \negmedspace : {1End}(x) &   \to   & Set     \\
                        \phi\in 2Hom(A, B)  & \mapsto & \phi_*, \\
  \end{eqnarray*}
  where
  $$
    \phi_*\negmedspace: \mathbb{T}r(A)\rightarrow\mathbb{T}r(B)
  $$
  is given by composition with $\phi$.
  A priori, $\mathbb{T}r$ is set valued, but if we assume $\mC$ to be enriched
  over a category $\mathcal{A}$ (cf.\ Definition \ref{enriched-Def}),
  $\mathbb{T}r$ takes values in $\mathcal{A}$.
  We will often assume that $\mC$ is $k$-linear for a fixed field $k$.
\end{Rem}

\subsection {Examples of the categorical trace}
Our first example is the motivational example mentioned above.
\begin{Exa}[$2$-vector spaces]
  Let $\mathcal {C}=\tvect$ and $x=[n]$. Then $A$ is an $n\times n$ matrix
  $A= (A_{ij})$, where the $A_{ij}$ are vector spaces. In this case,
  $$\mathbb{T}r(A) = \bigoplus_{i=1}^n A_{ii}.$$
\end{Exa}
\begin{Exa}[Categories] Let
  $\mathcal{C}=\mathcal{C}at$ and $x=\mathcal {V}$ a category, so $A:
  \mathcal{V}
  \rightarrow\mathcal{V}$ is an endofunctor. Then $\mathbb{T}r(A) =
  NT(\id_{\mathcal{V}}, A)$
is the set of natural transformations from the identity functor to $A$.
\end{Exa}
\begin{Exa}[Bimodules]
  Let $\mathcal{C}=\mathcal{DB}im$, so that $x=R$ is a ring, and $A=M$
  is an
  $R$-bimodule. Then
  $$\mathbb{T}r^\bullet(A) = Ext^\bullet_{R\otimes R^{op}}(R, M)$$
  is the Hochschild cohomology of $R$ with coefficients in $M$, see \cite{Loday}.
\end{Exa}
\begin{Exa}[Varieties]\label{Varieties-Exa} Let 
  $\mathcal{C}=\mathcal{V}ar_k$, $x=X$ be  a variety, and  $A=\mathcal{K}$ be a complex
  of coherent sheaves on $X\times X$. Then
  $$\mathbb{T}r^\bullet(A) = \mathbb{H}^\bullet(X,i^!(\mathcal {K})).$$
  here $i: X\rightarrow X\times X$ is the diagonal embedding, $i^!$ is
  the right adjoint of $i_*$, 
  and $\mathbb{H}$
  is the hypercohomology. In particular, if $\mathcal{K}$ is a vector bundle
  on $X\times X$ situated in degree 0, then
  $$\mathbb{T}r^\bullet(A) = H^\bullet(X, \mathcal{K}|_\Delta)$$
  is the cohomology of the restriction of $\mathcal{K}$ to the diagonal.
\end{Exa}

\subsection{The center of an object} The set $\mathbb{T}r({1}_x)$ will be called the
center of $x$ and denoted $Z(x)$. It is closed under both compositions $\circ_0$
and $\circ_1$. The following fact is well known \cite{Maclane}.
\begin{Prop}
\label{COMMUTATIVITY-OF-CENTER} The operations $\circ_0$
and $\circ_1$ on $Z(x)$ coincide and make it into a commutative monoid.
\end{Prop}

Thus, if $\mathcal {C}$ is pre-additive, then $Z(x)$ is a commutative ring and for
each $A: x\rightarrow x$ the group $\mathbb{T}r (A)$ is a $Z(x)$-module.

\subsection{Examples} (a) If $\mathcal{C} = \mathcal{A}b$ and $x=\mathcal{V}$
is an abelian category, then $Z(\mathcal{V})$, i.e., the ring of natural transformations
from the identity functor to itself is known as the Bernstein center of $\mathcal{V}$,
see \cite{Deligne}

\vskip .3cm
\noindent
(b) If $\mathcal{C}=\Pi(X)$ is the Poincare 2-category of a CW-complex $X$, then
$Z(x)= \pi_2(X,x)$ is the second homotopy group. Proposition
\ref{COMMUTATIVITY-OF-CENTER}  is the categorical analog of the
commutativity of
$\pi_2$.

\subsection{Conjugation invariance of the categorical trace}
\label{conj-inv}
In this section we assume for simplicity that the \tcat $\mC$ is strict.
Recall \cite{Maclane} that a 1-morphism $B: y\rightarrow x$ is called
an equivalence
if there exist a \omor $C: x\rightarrow y$ called {\em quasi-inverse}
and 2-isomorphisms $u: {1}_x\Rightarrow BC$,
$v: {1}_y\Rightarrow CB$.
For any object $x$, the $1$-morphism $B=1_x$ is an equivalence with
$C=1_x$ and $u,v$ the isomorphisms from \ref{2cat-Sec} (4). If $B$
and $B'$ are composable
$1$-morphisms, which are equivalences with quasi-inverses $(C,u,v)$ and
$(C',u',v')$ respectively, then $B''\circ_0 B$ is an equivalence with
quasi-inverse
\begin{equation}\label{composition-Eqn}
  \(C\circ_0 C',(B'\circ_0u\circ_0 C')\circ_1 u',(C\circ_0
  v'\circ_0 B)\circ_1 v\).
\end{equation}
\begin{Prop}[Conjugation invariance] \label{INVARIANCE-OF-TRACE}
(a) Let $A: x\rightarrow x$ be a 1-endomorphism and $B: x\rightarrow y$
 an equivalence with quasi-inverse $C$.
Then the rule
$$(\phi: {1}_x\Rightarrow A) \mapsto (B\circ_0 \phi\circ_0 C)\circ_1 u$$
defines a bijection of sets
$$\map{\psi(B, C, u,v)}{\mathbb{T}r(A)}{\mathbb{T}r(BAC)}.$$
By abuse of notation, we will write $\psi(B)$ when $C$, $u$ and $v$
are clear from the context.

\noindent
(b) Assume
that $B$ and $B'$ are 1-endomorphisms of $x$ and that both of them are
equivalences.
Then we have
$$\psi(B'\circ_0 B) = \psi(B')\circ\psi(B).$$

\noindent
(c) We have $\psi(1_x) = \id$.
\end{Prop}
\begin{proof}  (a) To explain the formula, note that we can view $u$
as a 2-morphism ${1}_y\Rightarrow BC= B\circ {1}_x\circ C$, while
$$B\circ_0\phi\circ_0 C: B\circ {1}_x\circ C\Rightarrow B\circ A\circ C.$$
Since $u$ is a 2-isomorphism, composing with $u$ is a bijection.

\vskip .2cm

\noindent
(b) This follows from the definition of $\psi(B)$ together with
(\ref{composition-Eqn}).

\vskip .2cm
\noindent
(c) Obvious.
\end{proof}

%
\begin{Prop}\label{TRACE-OF-DIRECT-SUM}
 Let $\mathcal{C}$ be an additive 2-category and $A, A': x\rightarrow x$
be two 1-morphisms. Then
$$\mathbb{T}r(A\oplus A') = \mathbb{T}r(A)\oplus \mathbb{T}r(A').$$
\end{Prop}

This is an immediate consequence of the fact that $\mathcal{H} =
Hom_{\mathcal{C}}(x,x)$ is an 
 additive category, and that therefore

$$Hom_{\mathcal{H}}({1}_x, A\oplus A') = Hom_{\mathcal{H}}({1}_x, A)
\oplus Hom_{\mathcal{H}}({1}_x, A').$$

\subsection {The joint trace} In the situation of Proposition \ref{conj-inv}
assume that $A$ and $B$ {\em commute}, i.e., that we are given a 
2-isomorphism
$$\eta: B\circ A\Rightarrow A\circ B.$$
 Then we have a map
$$B_*: \mathbb{T}r(A)\to \mathbb{T}r(A),$$
defined as the composition
$$\mathbb{T}r(A)\buildrel \psi(B)\over\longrightarrow \mathbb{T}r(BAC)
\buildrel \mathbb{T}r(\eta\circ_0 1)\over\longrightarrow \mathbb{T}r(ABC)
\buildrel \mathbb{T}r(1\circ_0 u^{-1})\over\longrightarrow\mathbb{T}r(A).$$
Assume now that the 2-category $\mathcal{C}$ is $k$-linear for a field $k$.
Then $\mathbb{T}r(A)$ is a $k$-vector space, and $B_*$ is a linear operator.
Let us further assume that $\mathbb{T}r(A)$ is finite-dimensional.
Then we define the {\em joint trace} of $A$ and $B$ to be the following
element of $k$:
$$\tau(A, B) = \tr\{B_*: \mathbb{T}r(A)\to \mathbb{T}r(A)\}.$$
It depends on  the choice of the commutativity isomorphism $\eta$, as well
as on the equivalence data for $\Phi$. 

\section{2-representations and their characters}
\subsection{$2$-representations}
Let $G$ be a group. We view
$G$ as a $2$-category with one object, $\pt$, the set of 1-morphisms
$Hom(\pt, \pt)=G$ and all the  2-morphisms being the identities of the above
1-morphisms.

\begin{Def}\label{2-rep-Def}
  Let $\mC$ be a $2$-category. A 2-representation of $G$ in $\mC$ is a
  strong $2$-functor from $G$ to $\mC$. More explicitly, this is a system
of the following data:

  \begin{enumerate}
    \renewcommand{\labelenumi}{(\alph{enumi})}
    \item an object $V$ of $\mC$,
    \item for each element $g\in G$, a $1$-automorphism $\rho(g)$ of
      $V$,
    \item for any pair of elements $(g,h)$ of $G$ a $2$-isomorphism
      $$
        \tisomap{\phi_{g,h}}{(\rho(g)\circ\rho(h))}{\rho(gh)},
      $$
    \item and a $2$-isomorphism
      $$
        \tisomap{\phi_1}{\rho(1)}{\id_c},
      $$
  \end{enumerate}
  such that the following conditions hold
  \begin{enumerate}
    \renewcommand{\labelenumi}{(\alph{enumi})}
    \addtocounter{enumi}{4}
    \item for any $g,h,k\in G$ we have
      $$
        \phi_{(gh,k)}(\phi_{g,h}\circ\rho(k)) = \phi_{(g,hk)}(\rho(g)\circ\phi_{h,k})
      $$
      (associativity); we also write $\phi_{g,h,k}$,
    \item we have
      $$
        \phi_{1,g} = \phi_1\circ\rho(g) \quad\text{and}\quad
        \phi_{g,1} = \rho(g)\circ\phi_1.
      $$
  \end{enumerate}
\end{Def}

Note that this definition is the special case of the concept of
a representation of a  $2$-group
as defined by  Elgueta \cite[Def.4.1]{Elgueta:2groups}. 
This case corresponds to the 2-group being discrete, i.e., being reduced to an
ordinary group. Compare also \cite[\S 0]{Del}. 
If $\mathcal D$ and $\mathcal C$ are $2$-categories, then  
strong $2$-functors from $\mathcal D$ to $\mathcal C$ form a
$2$-category $\frak{Hom}(\mathcal{D}, \mathcal{C})$, see
 \cite[Def. I.1.9]{Hakim} for strict 2-categories or  \cite{Benabou}
for the general case. 
In particular, 
$2$-representations of $G$ in $\mathcal{C}$ form a 2-category. We
will denote it by $2Rep_{\mathcal{C}}(G)$.
We understand that the implications of this fact will be spelled out in
detail in \cite{Bartlett:Willerton}. 

\subsection{The category of equivariant objects} 
Consider the particular case of
Definition \ref{2-rep-Def} when $\mathcal{C}=\mathcal{C}at$ is the 2-category of
(small) categories. Then a 2-representation of $G$ in $\mathcal{C}$ is the same as an
action of $G$ on a category $\mathcal{V}$. In other words, each
$$\rho(g): \mathcal{V}\to\mathcal{V}$$
is a functor and each $\phi_{g,h}$ is a natural transformation.
We will call a category with  a $G$-action a {\em categorical representation}
of $G$ and will denote by $2Rep (G)  = 2Rep_{\mathcal{C}at}(G)$
the 2-category formed by  categorical representations.

In this section, we formulate the categorical analogue 
the concept of the  subspace of $G$-invariants of a representation.
\begin{Def}
  Fix a category $\underline 1$ with one object and one
  morphism. 
  The {\em trivial $2$-representation
    of $G$} is given by the unique action of $G$ on $\underline 1$.
  We will also denote it by $\underline 1$.
  Let $\rho$ be an action of $G$ on $\mV$. We define the {\em
  category of  $G$-equivariant objects in $\mV$}, denoted $\mV^G$, to be the category
  of $G$-functors  from $\underline 1$ to $\rho$:
$$\mV^G = Hom_{2Rep(G)} (\underline{1}, \rho).$$
\end{Def}
This definition spells out to the following. 
An object of $\mathcal{V}^G$ consists of an object $X\in\ob(\mV)$ and a system
$$\( \map{ \epsilon_g}{X}{\rho(g)(X)},\medspace g\in G\),$$
where $\epsilon_g$ are isomorphisms satisfying the following
compatibility conditions: First, it is required that for $g=1$ we have
$$\epsilon_1 = \phi_{1,X}^{-1}: X\mapsto \rho(1)(X).$$
Second, it is required that for any $g,h\in G$ the diagram
$$\xymatrix{
X\ar[r]^{\epsilon_g}\ar[d]_{\epsilon_{gh}} &
{\rho(g)(X)} \ar[d]^{\rho(g)(\epsilon_h)} \\
{\rho(gh)(X)} &
{\rho(g)(\rho(h)(X))}\ar[l]^{\phi_{g,h,X}}
}$$
is commutative.

\begin{Exa}\label{trivial} Let $\mathcal W$ be a category.
 Define the trivial action of $G$ on $\mathcal W$ by taking
  all $\rho(g)$ and $\phi_{g,h}$ to be the identities.
Then a $G$-equivariant object in $\mathcal{W}$ is the same as a representation of $G$
in $\mathcal W$, i.e., an object $X\in\mathcal{W}$ and a homomorphism
$G\to \on{Aut}_{\mathcal {W}}(X)$. 
\end{Exa}

\begin{Prop}
  Let $\mathcal W$ be a category equipped with trivial $G$-action as in Example \ref{trivial}.
 Then we have an equivalence of categories
$$Hom_{2Rep(G)}(\mathcal{W}, \mathcal {V}) \simeq Hom_{\mathcal{C}at}(\mathcal{W}, \mathcal{V}^G).$$
In particular, taking $\mathcal{W}=\on{pt}$ (the category with one 
object and one morphism), we get
$$\mathcal{V}^G \simeq Hom_{2Rep(G)}(\on{pt}, \mathcal{V}).$$
\end{Prop}

In plain words, this means that  any $G$-functor from $\mathcal W$
 to $\mathcal V$ factors through the forgetful functor
$$i_{\mathcal V}: \mathcal V^G \to \mV.$$
In 2-categorical terms, this can be formulated by saying that the 2-functor
$$I_G: 2Rep(G)\to \mathcal{C}at, \quad \mathcal{V}\mapsto \mathcal{V}^G,$$
is right 2-adjoint (in the sense of \cite[Def. I.1.10]{Hakim}),
to the 2-functor $\mathcal{C}at\to 2Rep(G)$ associating
to any $\mathcal{W}$ the same category $\mathcal{W}$  with trivial $G$-action.

\begin{Pf}{} This follows at once from the definition of the $Hom$-categories
in $2Rep(G)$, which are particular cases of $Hom$-categories in 2-categories
of 2-functors, see {\em loc. cit.} Def. I.1.9. Indeed,
denote by $\widetilde{\rho}$ the trivial action of $G$ on $\mathcal{W}$. 
Then a $G$-functor
$F: \mathcal{W}\to\mathcal{V}$ gives, for each object $X\in \mathcal{W}$,
an object $F(X)\in\mathcal{W}$ together with isomorphisms
$$u_{g,X}: F(\widetilde{\rho}(g)(X))\to \rho(g)(F(X)),$$ 
satisfying the compatibility condition for each pair $g,h\in G$.
Since  $\widetilde{\rho}(g)(X)=X$, the system formed by $F(X)$ and
the $u_{g,X}$ gives an equivariant object of $\mathcal{V}$. We leave 
further details to the reader. \end{Pf}

\begin{Rem} The concept of the category of equivariant objects relates
our approach to 2-representations with a different approach  due to Ostrik
 \cite{Ostrik}. If $k$ is a field, then finite-dimensional linear representations
of $G$ over $k$ form a monoidal  category 
 $(\mathcal{R}ep(G),  \otimes)$
with respect to the usual tensor product. In {\em loc. cit.} it was proposed
to study module categories over $\mathcal{R}ep(G)$. In our situation,
given a $G$-action on a $k$-linear additive category $\mathcal{V}$, the
category $\mathcal{V}^G$ is naturally a module category over $\mathcal{R}ep(G)$.
In other words, the tensor product of a $G$-representation and a $G$-equivariant
object is again a $G$-equivariant object.  It seems that in general, the
passage from a $G$-category $\mV$ to the $\mathcal{R}ep(G)$-module category
$\mV^G$ leads to some loss of information. However, in some
particular cases, the two approaches are equivalent, see Remark \ref{Ost}
below.

\end{Rem}

\subsection{Characters of  $2$-representations}

We are now ready to define the categorical character of a
$2$-representation. To motivate the dicussion of this section, we
start with a
reminder of classical character theory.
\subsubsection{Group characters and class functions}
We fix a field $k$ of characteristic $0$
containing all roots of unity.
 Let $G$ be a group.
Recall that a function $f\negmedspace : G\rightarrow k$ is called a
{\em class function} if
it is invariant under conjugation.
\begin{Not}
We denote by $\on{Cl}(G; k)$
the ring of class functions on $G$. As before, let
 $\mathcal{R}ep(G)$ be the
category of finite-dimensional representations of $G$ over
$k$. We write $R(G)$
for its Grothendieck ring $K(\mathcal{R}ep(G))$.
\end{Not}
If $\rho: G\rightarrow Aut(V)$ is a representation, then its
character
\begin{eqnarray*}
  \chi_V\negmedspace : G &   \to    &  k \\
                       g &  \mapsto &  \tr(\rho(g))
\end{eqnarray*}
is  a class function. The following is well known \cite{Serre:representations}.
\begin{Prop} \label{CLASS-FUNCTIONS}
If $G$ is finite, then the correspondence $V\mapsto \chi_V$
induces an isomorphism of rings
$$R(G)\otimes k \rightarrow \on{Cl}(G; k).$$
\end{Prop}
\subsubsection{The categorical character}
The classical definitions discussed in the previous section suggest
the following analogues for $2$-representations:
\begin{Def}
  Let $\rho$ be a \trep of $G$.
  We define the \tchar of $\rho$ to be the assignment
  $$
    g \mapsto \ttr{\rho(g)}.
  $$
\end{Def}
We now discuss the sense in which the \tchar is a class function.
First we recall the definition of the inertia groupoid of $G$:
\begin{Def}\label{inertia-Def}
  Let $G$ be a group. The {\em inertia
    groupoid}
  $\LG$ of $G$ is the category that has as objects the
  elements of
  $G$ and
$$Hom_{\LG}(u,v) = \bigl\{ g\in G: \ v=gug^{-1}\bigr\}.$$

\end{Def}
\begin{Prop}\label{class-fcn-Prop}
  Let $\mC$ be a \tcat
  and let $\rho$ be a 2- representation of $G$ in $\mC$.
  Then the \tchar of $\rho$ is a functor from the inertia
  groupoid $\Lambda(G)$ to the category of sets:
  $$
    \map{\ttr{\rho}}{\Lambda(G)}{\Set}.
  $$
  If $\mC$ is enriched over a category $\mathcal A$, then this functor takes
  values in $\mathcal A$. In other words, for any $f,g\in G$ there is an isomorphism
 $$\map{\psi(g)=\psi_f(g)}{\ttr{\rho(f)}}{\ttr{\rho(gfg^{-1})}},$$
and these isomorphisms satisfy
 \begin{enumerate}
    \renewcommand{\labelenumi}{(\alph{enumi})}
    \item $\psi(gh) = \psi(g)\circ\psi(h)$ and
    \item $\psi(1) = \id_{\rho(f)}$.
  \end{enumerate}
\end{Prop}
\begin{proof}{}
  Pick $f,g\in G$ and write ${A}$ for $\rho(f)$, $B$ for $\rho(g)$, $C$
  for $\rho(g^{-1})$ and define $$\map{u}{1_c}{{BC}}$$ as the composite
  of maps from Definition \ref{2-rep-Def}: 
  $$u := \phi^{-1}_{g,g^{-1}}\phi^{-1}_1.$$ With this notation,
  Proposition \ref{INVARIANCE-OF-TRACE} (a) implies the existence of
  an isomorphism
  $$
    \isomap{\psi'}{\ttr{\rho(f)}}{\ttr{\rho(g)\rho(f)\rho(g^{-1})}}.
  $$
  Composed with $\ttr{\phi_{g,f,g^{-1}}}$ this gives the desired map
  $\psi(g)$. Properties (a) and (b) of $\psi(g)$ follow from
  Proposition \ref{INVARIANCE-OF-TRACE} (b) and (c).
\end{proof}

\begin{Rem}\label{sheaf} By regarding $G$ as a discrete topological space,
we can consider the correspondence
$g\mapsto \ttr{\rho(g)}$ as a sheaf of sets on $G$. 
If $G$ is a topological or algebraic
group,  there are natural situations when
 $\ttr{\rho}$ is a sheaf on $G$
in the corresponding stronger sense,  equivariant under
conjugation, see Subsection \ref{Contructible-sheaves-Sec} below
for an example.  
\end{Rem}

\begin{Def}
  If $\rho$ is a \trep in a $k$-linear \tcat
with finite-dimensional $\thom(\phi,\psi)$, we define the \tchar of
  $\rho$ to be the function $\chi_\rho$ on pairs of commuting
  elements
given by the joint trace of $\rho(g)$ and $\rho(h)$: 
$$
 \chi_\rho (g,h)=\tau(\rho(g), \rho(h))  =  \on{Trace}\bigl\{ \psi(h):\ttr{\rho(g)}
\to \ttr{\rho(g)}\bigr\}.
$$

\end{Def}

Note that

\begin{equation}\label{2cl}
\chi_\rho(s^{-1}gs, s^{-1}hs) = \chi_\rho(g,h).
\end{equation}

This can be formulated as follows.

\begin{Def} Let $G$ be a group and $R$ be a commutative ring.
A 2-class function on $G$ with values in $R$ is a function
$\chi(g,h)$ defined on pairs of commuting elements of $G$ and
invariant under simultaneous conjugation, as in
\eqref{2cl}. The ring of such functions will be denoted $2\operatorname{Cl} (G; R)$.
\end{Def}
Thus the \tchar is a 2-class function.

\section{Examples}

\subsection{$1$-dimensional  $ 2$-representations}\label{1-dim-Sec}

Let $k$ be a field and
$$c: G\times G\to k^*$$ be a 2-cocycle,
i.e., it satisfies the identity
$$c(g_1 g_2, g_3) c(g_1, g_2) = c(g_1, g_2 g_3) c(g_2, g_3).$$
We then have an action $\rho=\rho_c$ of $G$ on $Vect_k$.
By definition, for $g\in G$ the functor $\rho(g): Vect_k\to Vect_k$
is the identity, while
$$\phi_{g,h}\negmedspace : \id=\rho(g)\circ\rho(h) \Rightarrow
\rho(gh)=\id$$
is the multiplication with $c(g,h)$, and $\phi_1$ is the
multiplication by $c(1,1)$.
The cocycle condition for $c$ is equivalent to Condition (e) of
Definition \ref{2-rep-Def}, while Condition (f) follows because
$$
  c(1,1g)\cdot c(1,g)=c(1,1)\cdot c(1,g)
$$
implies that $$c(1,g)=c(1,1)$$
and similarly $$c(g,1)=c(1,1).$$
Cohomologous cocycles define equivalent $2$-representations, and it is
easy to see that $H^2(G,k^*)$ is identified with the set of
$G$-actions on $\on{Vect}_k$ modulo equivalence. Compare \cite{Kapranov}.

\vskip .2cm

We now find the categorical character and the $2$-character of
$\rho_c$. First of all, the functor $\rho_c(g)$ being the identity,
$$\ttr{\rho_c(g)}=k.$$
Next, the equivariant structure on $\ttr{\rho_c}$ was defined in the
proofs of Propositions
\ref{INVARIANCE-OF-TRACE} (a) and \ref{class-fcn-Prop} to be the
composition
$$
  \ttr{\rho_c(f)}\buildrel{\tilde u}\over{\longrightarrow} \ttr
  {\rho_c(g)\rho_c(f)\rho_c(g^{-1})} \longrightarrow
  \ttr{\rho_c(gfg^{-1})}.
$$
Here $\tilde u$ is induced by the $1$-composition with
$$u\negmedspace :1_{\on{Vect}_k}\Rightarrow BC,$$
where $B=\rho_c(g)$ and $C=\rho_c(g^{-1})$.
In our case,
$$u=c(g,g^{-1})^{-1}$$
(multiplication with a scalar).
The second map is induced by
$$\phi_{g,f,g^{-1}}=c(g,f)c(gf,g^{-1}),$$
see Definition \ref{2-rep-Def} (e). As a result we have
\begin{Prop}
  For any two commuting elements $f,g\in G$, we have
  $$
    \chi_{\rho_c}(f,g) = c(g,f)c(gf,g^{-1})c(g,g^{-1})^{-1}.
  $$
\end{Prop}
Notice also the following fact which extends Example \ref{trivial}. 
\begin{Prop}[{compare \cite{Elgueta:2groups}}]
  Let $\rho_c$ be the one-dimensional \trep of $G$ on $\mV=\on{Vect}_k$
  corresponding to $c$. Then objects of $\mV^G$ are the same as projective
  representations of $G$ with central charge $c$, i.e., pairs
$(V, \varphi: G\to\on{Aut}(V))$, where $V$ is a $k$-vector space, and 
 $\varphi$ is a map satisfying
$$\varphi(gh) = \varphi(g)\varphi(h)\cdot c(g,h).$$
\end{Prop}

\subsection{Representations on $2$-vector spaces}
\treps $\rho_c$ from Section \ref{1-dim-Sec}
can be viewed as acting on the $1$-dimensional $2$-vector space
$[1]$. More precisely, let
$$
  1\to k^*\to\widetilde {G}\stackrel\pi\to G\to 1
$$
be the central extension corresponding to the cocycle $c$. For every
$g\in G$ the set $\pi\inv(g)$ is the a $k^*$-torsor, and therefore
$$
  L_g := \pi\inv(g)\cup\{0\}
$$
is a $1$-dimensional $k$-vector space. The group structure on $\widetilde
{G}$ induces isomorphisms
$$
  L_g\tensor_k L_h\to L_{gh}.
$$
It follows that associating to $g\in G$ the $2$-matrix
$\| L_g\|$ of size $1\times 1$ gives a \trep of $G$ on $[1]\in\ob(\tvect)$

More generally, a \trep $\rho$ of $G$ on $[n]$ consists of the
following data: for
each $g\in G$, a quasi-invertible $2$-matrix
$\rho(g)=\|\rho(g)_{ij}\|$ of size $n\times n$, with each $\rho(g)_{ij}$ being a
$k$-vector space, plus the data $\phi_{g,h}$ as in Definition \ref{2-rep-Def}.
\begin{Lem}\label{perm-Lem}
  A $2$-matrix $A=\|A_{ij}\|$ of size $n\times n$ is quasi-invertible
  if and only if there is a permutation $\sigma\in\Sn$ such that
  $A_{ij}=0$ for $i\neq\sigma(j)$ and $\dim(A_{i,\sigma (i)})=1$.
\end{Lem}
It follows that an $n$-dimensional \trep of $G$ defines a homomorphism
$G\to\Sn$ plus some cocycle data. This is naturally explained in the
context of induced 2-representations, cf.\ Section \ref{induced-Sec} below.
\begin{Rem}
  Lemma \ref{perm-Lem} suggests that the theory becomes richer if
  one
  works with generalizations of $2Vect$ that have more interesting
  quasi-invertible 1-morphisms. One such generalization was defined in
  \cite{Elgueta:generalized.linear.2.groups}.
\end{Rem}
\subsection{Constructible sheaves}\label{Contructible-sheaves-Sec}
Let $X$ be a real analytic manifold acted upon by  $G$. We view each $g\in
G$ as  map  $\map gXX$. Let $\mathcal V$ be the category
$\mathcal D^bConstr(X)$, see Section \ref{Examples} (d). The base field
$k$ will be taken to be the field $\CC$ of complex numbers, for simplicity.
 We have
then an action $\rho$ of $G$ on $\mathcal V$ given by
$$
  \rho(g)(\mathcal F) = (g\inv)^*(\mathcal F)
$$
(inverse image under $g^{-1}$).  
As in Examples \ref{Examples} (c), (d), it is more practicable to lift
this action on a category to an action on an object $X$ of the 2-category
${\Bbb R}{\mathcal A}n_{\CC}$. Namely, for  $g\in G$ we denote by
$\Gamma(g)\subset X\times X$ its graph and associate to $g$ the constructible
sheaf $\mathcal{K}_g = \underline{\CC}_{\Gamma(g)}$, the constant sheaf on $\Gamma(g)$.
Note that $\rho(g) = F_{\mathcal{K}_g}$ is the functor associated to $\mathcal{K}_g$.
It is clear that the correspondence $g\mapsto \mathcal{K}_g$ gives an action
of $G$ on the object $X$ which we denote
  $\widetilde{\rho}$.

\begin{Prop}\label{sheaves}
  Assume that $X$ is  oriented.
  Then the categorical character of $\widetilde{\rho}$ is found as follows:
  $$
    \mathbb Tr^\bullet(\widetilde{\rho}(g)) = H^{\bullet+\on{codim} X^g}(X^g,\CC),
  $$
  where $X^g\sub X$ is the fixed point locus of $g$.
\end{Prop}
\begin{proof}
As in Example \ref{Varieties-Exa}, we denote by $i: X\to X\times X$ the
diagonal embedding and we
  have
  $$
  \ttr{\widetilde{\rho}(g)}=\ttr{F_{\mathcal K_g}} = \mathbb H^\bullet(X,i^!
(\mathcal K_g))
  $$
  and
  $$i^!(\mathcal{K}_g) = \underline{R\Gamma}_{\Delta}(\mathcal{K}_g) =
 \underline{ \CC}_{X^g}[\on{codim} (X^g)].$$
  Here $\Delta = i(X)$ is the diagonal in $X\times X$.
\end{proof}

Assume further that $X$ is compact, and $G$ is a Lie group acting smoothly
  on $X$. Then Proposition \ref{sheaves} can be sheafified as follows. 
Let 
$$Y=\bigl\{ (g,x)\in G\times X| \, gx=x\}$$ 
be the ``universal'' fixed point space. 
We have the natural embedding and projection
\begin{equation}\label{lus}
G\times X\buildrel \eta\over\longleftarrow Y \buildrel p\over\longrightarrow G.
\end{equation}
Further, the group $G$ acts on the left on $G\times X$, preserving $Y$, by the formula
$$g_0 (g,x) = (g_0 g g_0^{-1}, g_0x).$$
Thus $\eta$ is $G$-equivariant, and so is $p$, if we consider the action of $G$
on itself by conjugations. 
Thus the constructible complex of sheaves
\begin{equation}\label{T}
\frak{T}^\bullet  = Rp_*   \underline{\CC}_{Y}
\end{equation}
on $G$ is conjugation equivariant. 
It can be seen as a sheaf-theoretical version of the
 categorical trace. Indeed, for any $g\in G$ the complex $\frak{T}^\bullet_g$, the
stalk of $\frak{T}^\bullet$ at any $g\in G$ 
has cohomology
\begin{equation}\label{stalk}
H^i(\frak{T}^\bullet_g) = H^i(X^g, \CC) = 
 \mathbb Tr^{\bullet - \on{codim}(X^g)} (\widetilde{\rho}(g)).
\end{equation}

\begin{Exa}[Character sheaves] We specialize to the case when $G$ is a complex
semisimple algebraic group and $X$  is the flag variety of $G$,
i.e., the space of all Borel subgroups in $G$
with $G$-action by conjugation. In this case a class of
conjugation equivariant complexes on $G$ was constructed by G. Lusztig in the
framework of his theory of character sheaves \cite{Lusztig}. Lusztig's complexes
are grouped into ``series'' labeled by an element $w$ of the Weyl group.
We consider here the ``principal'' series corresponding to $w=1$. Complexes
of this series are defined in terms of the diagram (3) and can be interpreted
as categorical traces of certain 2-representations of $G$, via a twisted version of
Proposition \ref{sheaves} and the equality \eqref{stalk}.

To be precise, recall  that  all Borel subgroups $B\subset G$
are  conjugate and the normalizer of any $B$ is $B$ itself. Therefore,  
 the abelianizations
$B/[B,B]$ for different $B$  are canonically identified with each other. Equivalently, we
can say that they are all identified with a fixed group $T$ (the ``abstract''
maximal torus, cf. \cite{Chriss:Ginzburg}, p. 137). Since in our case
$$Y= \{ (g, B)\in G\times X: g\in B\},$$
we get  a projection $q: Y\to T$ taking $(g,B)$ to the image of $g$
in the abelianization of $B$.  Given a 1-dimensional
local system $\mathcal L$ on $T$, the sheaf $q^*{\mathcal L}$ on $Y$ is $G$-equivariant,
and hence the constructible complex 
\begin{equation}\label{TL}
\frak{T}^\bullet (\mathcal {L}) =  Rp_* q^* {\mathcal L}
\end{equation}
on $G$ is conjugation invariant. The complex from \eqref{T} corresponds to ${\mathcal L}=\underline{\CC}_T$. 
 Lusztig's character sheaves (corresponding to
$w=1$) are direct summands of $\frak{T}^\bullet (\mathcal L)$ in the
derived category.

This can be interpreted as follows. Let $\pi: Z\to X$ be the basic affine space of
$G$. It is a principal $T$-bundle on $X$ isomorphic to the product of the  $\CC^*$-bundles
corresponding to fundamental weights,  see,  e.g., \cite{Braverman:Polishchuk}, Sect. 2.3.2.
The fibers $\pi^{-1}(B), B\in X$,  are identified with $T$ uniquely up to a
group translation in $T$. Although the local system $\mathcal{L}$ on $T$ is
not $T$-equivariant (unless it is trivial),
 each translation takes it into an isomorphic sheaf. 
Therefore
 it makes sense to speak about local systems on $\pi^{-1}(B)$
isomorphic to  $\mathcal{L}$ on $T$.

Call an $\mathcal{L}$-{\em twisted sheaf} on $X$ a sheaf on $Z$ whose restriction
on each fiber of $\pi$ has the form $\mathcal{L}'\otimes V$ where $\mathcal{L}'$
is a local system isomorphic to $\mathcal{L}$ on $T$, and $V$ is a $\CC$-vector space.
Let 
$$\mathcal{V}(\mathcal{L}) = D^b_{\mathcal{ L}} Constr (X)$$
be the derived category of bounded complexes of $\mathcal L$-twisted sheaves on $X$ with
constructible cohomology. 
We have a natural action $\rho_{\mathcal L}$ of  $G$  on $\mathcal{V}(\mathcal{L})$, and the
stalk of $\frak{T}^\bullet({\mathcal L})$ at $g$ can be related to the categorical trace
of $\rho_{\mathcal L}(g)$  similarly to \eqref{stalk}. As before, to make this
precise, we need to lift $\rho_{\mathcal L}$
to an action $\widetilde{\rho}_{\mathcal L}$  by ``kernels''  and restrict the
 kernels to the diagonal.
Compare \cite{Braverman:Polishchuk}, Sect. 3.3-4.

\vskip 1cm

\end{Exa}

\section{Representations of finite groupoids}

\subsection{Reminder on semisimplicity}

Recall that a groupoid is a category with all morphisms invertible.
A groupoid is called finite if it has finitely many objects and morphisms.
As before, we fix  a field  $k$ of characteristic $0$
containing all roots of unity. We denote by $\Vect_k$ the category
of finite-dimensional $k$-vector spaces. 

\begin{Def}
  Let $G$ be a finite groupoid.
  A  (finite-dimensional) {\em representation over $k$} of $G$ is a functor from $G$ to
  $\Vect_k$. A morphism between two $G$-representations is a
  natural transformation between them.
  We denote the category of $G$-representations over $k$ by
  $\Repk G$.
  Object-wise direct sum and  tensor product make $\Repk G$ into a
  bimonoidal category, so that its Grothendieck group $R(G)$
(with respect to the direct sum)  is
  a ring, the {\em representation ring} of $G$.
\end{Def}

\begin{Def}
  Let $\map{\alpha}{H}{G}$ be a map of finite groupoids. Then
  precomposition with $\alpha$ defines a functor
  $$
    \map{\res{}\alpha}{\Repk G}{\Repk H}.$$
  If $H$ is a subgroupoid of $G$ and $\alpha$ is its inclusion we also
  denote $\res{}{\alpha}$ by $\res GH$.
\end{Def}
\begin{Def}
  Let $G$ be a finite groupoid. The groupoid algebra $k[G]$ has
  as underlying $k$-vector space
  the vector space with one basis-element $e_g$ for each morphism $g$ of
  $G$. The algebra structure is given by
  $$e_g\cdot e_h =
  \begin{cases}
    e_{gh} & \text{if $g$ and $h$ are composable}\\
    0& \text{else}
  \end{cases}
  $$
\end{Def}

The categories of $G$-representations over $k$ and of $k[G]$-modules
are then equivalent.
\begin{Prop}
  If $\map\alpha HG$ is an equivalence of groupoids, then
  $$\map{\res{}\alpha}{\Repk G}{\Repk H}$$ is an equivalence of categories.
\end{Prop}
\begin{Cor}
  If the groupoids $G$ and $H$ are equivalent, then their groupoid
  algebras are Morita equivalent.
\end{Cor}
\begin{Obs}\label{groups-Obs}
  Every finite groupoid is equivalent to a disjoint union of groups.
\end{Obs}
\begin{Pf}{}
  Let $G$ be a finite groupoid. Pick a representative for each
  isomorphism class of objects in $G$. Consider the inclusion of the
  disjoint union of the automorphism groups of these representatives in
  $G$. By construction, this inclusion is fully faithful and
  essentially surjective, so it is an equivalence of categories.
\end{Pf}
\begin{Cor}
  The groupoid algebra of a finite groupoid $G$ is semi-simple.
  Thus there is a unique decomposition  of representations of $G$ into
  irreducibles.
\end{Cor}

\begin{Def}\label{in-group}
 (a)  Let $G$ be a groupoid. The inertia groupoid $\LG$ of $G$
  has as objects, the automorphisms of $G$ (i.e., morphisms $u\in\on{Mor}(G)$
whose source and target coincide).
 For two such morphism $u,v$ there is one  morphism
in $\LG$  from $u$ to
  $v$ for every morphism $g$ of $G$ with $v=gug\inv$.

(b)  A class function on a groupoid $G$ is a
  function defined on isomorphism classes of objects
  of $\LG$.

(c)  Let $\rho$ be representation of a finite groupoid $G$.
  The character
  of $\rho$ is the $k$-valued class function on $G$ given by
  $$\chi_\rho([g]) = \tr(\rho(g)).$$
\end{Def}

 As before for the case of groups, we denote by
$\on{Cl}(G;k)$ the ring of class functions on $G$. 
\begin{Cor} 
 Sending a representation to its character is a ring map
  $\RG\to \Cl{G}.$
This map 
 becomes an isomorphism after
  tensoring with $k$.
\end{Cor}

\begin {Pf}{} The first statement  is obvious, the second one 
follows from Observation \ref{groups-Obs} and Proposition 
\ref{CLASS-FUNCTIONS}. 
\end{Pf}

\subsection{Induced representations of groupoids}

\begin{Def}\label{ind-Def}
  Let
  $
   \map{\alpha}HG
  $
  be a map of groupoids, and let $V$ be a representation of $H$.
  Viewing $V$ as a $k[H]$-module, we define the induced
  $G$-representation of $V$ by
  $$
    \ind\alpha{}(V) := k[G]\tensor_{k[H]}V.
  $$
  We will sometimes write $\ind HG$ for $\ind\alpha{}$, if the map is
  obvious.
\end{Def}
Note that $\ind\alpha{}$ is left adjoint to $\res{}\alpha$.

Let
$$\map{\alpha}HG$$
be faithful and essentially surjective,
and let $\rho$ be a
representation of $H$.
By Observation \ref{groups-Obs}, we may assume that $G$ and $H$ are
disjoint unions of groups. Then a
representation of $G$ can be described one group at a time, so we may
as well assume that $G$ is a single group and
$$\map{\alpha}{H_1\sqcup \dots \sqcup H_n}G$$
is given by a (nonempty) set of injective group maps
$\alpha_1,\dots,\alpha_n$.

  In this situation, the induced representation of $\rho$
  along $\alpha$ is isomorphic to
  $$
    \ind\alpha{}(\rho_1,\dots,\rho_n) := \bigoplus \ind{\alpha_i}{}\rho_i.
  $$

Note also that away from the essential image of $\alpha$, the induced
representation along $\alpha$ is zero.

\begin{Prop} Assume that  $\alpha$ is a faithful functor,
let $x$ be an object of $G$ and $g\in Hom_G(x,x)$. Then
the character of the induced representation evaluated at $g$ is given
by the formula
$$
  \chi_{\on{ind}}(x,g) = \sum_{y\in H_0} \frac{1}{|\on{orbit}_H(y)|}
  \frac{1}{|\Aut_H(y)|}
 \sum_{\stackrel{s\in G_1\mid sx = y}{sgs^{-1}\in\Aut_H(y)}}
  \chi(y,sgs^{-1}).
$$
Here $\on{orbit}_H(y)$ is the $H$-isomorphism class of $y$, and
 the second sum is over all morphisms $s$ of $G$ with source $x$
and target $y$ that conjugate $g$ into a morphism in the image of
$\alpha\at{\Aut_H(y)}$.
\end{Prop}
\begin{Pf}{}
  Let $\on{orbit}_G(x)$ denote the $G$-isomorphism class of $x$, and let
  $$\mathcal{R}=\{y_1,\dots,y_n\}$$
  be a system of representatives for the $H$-isomorphism classes
  mapping to $\on{orbit}_G(x)$ under $\alpha$.
  For each $j$, pick an $s_j$ with $s_jx=\alpha(y_j)$.
  Denote $\alpha\at{\Aut_H(y_j)}$ by $\alpha_j$.
  We have
  $$
    \chi_{\on{ind}} = \sum_{j=1}^n\chi_{\on{ind}\at{\alpha_j}}
  $$
  and
  $$\chi_{\ind{\alpha_j}{}}(x,g) = \chi_{\ind{\alpha_j}{}}(s_jx,s_jgs_j^{-1}).$$
  By the classical formula for the character of an induced
  representation of a group, we have
  \begin{eqnarray*}
  \chi_{\on{ind\at{\alpha_j}}}(s_jx,s_jgs_j^{-1})& =&
  \frac{1}{|\Aut_H(y_j)|}
  \sum_{\stackrel{t\in\Aut_G(y_j)}{ts_jgs_j^{-1}t^{-1}\in\Aut_H(y_j)}}
  \chi(y_j,ts_jgs_j^{-1}t^{-1})\\ & = &\frac{1}{|\Aut_H(y_j)|}
  \sum_{\stackrel{sx = y_j}{sgs^{-1}\in\Aut_H(y_j)}}\chi(y_j,sgs^{-1}).
  \end{eqnarray*}

  The first sum in the Proposition is over all
  objects $y$ of $H$. If $\alpha(y)$ is not isomorphic to $x$, then
  the second sum is empty. For the ones isomorphic to $x$, we are double
  counting: rather than just having one summand for the
  representative $y_j$, we have the same summand for every
  element in its $H$-orbit.
\end{Pf}
\begin{Cor}
  Consider an inclusion of groups $H\sub G$ and the induced map of
  inertia groupoids
  $\map\alpha{\Lambda(H)}\LG$.
  Let $\rho$ be a representation of $\Lambda(H)$. Then the character
  of the induced representation $\ind\alpha{}\rho$ evaluated at a pair
  of commuting elements of $G$ is given by the
  formula
  \begin{eqnarray*}
      \chi(g_1,g_2) & = &\sum_{h_1\in
      H}\frac{1}{|[h_1]_H|}\frac{1}{|C_H(h_1)|} \sum_{\stackrel{
        sg_1s^{-1}=h_1}{sg_2s^{-1}\in
        C_H(h_1)}}\chi(sg_1s^{-1},sg_2s^{-1}) \\
     & = & \frac{1}{|H|}\sum_{\stackrel{s(g_1,g_2)s^{-1}}{\in H\times H}}
     \chi(sg_1s^{-1},sg_2s^{-1}).
  \end{eqnarray*}
Here $[h_1]_H$ is the conjugacy class of $h_1$ in $H$ while $C_H(h_1)$
is the centralizer of $h_1$ in $H$.
\end{Cor}

\section{Induced 2-representations}\label{induced-Sec}

\subsection{Main definitions}

\begin{Def}
  Let $H\sub G$ be an inclusion of finite groups. Let
  $$\map{\rho}{H}{\on{Fun}(\mV, \mV)}$$
  be an action of $H$ on a category $\mV$.
  Let $\ind HG(\mV)$ be the category whose objects are maps
  $$\map fG\ob{\mV}$$
  together with an isomorphism for every $g\in G$ and $h\in H$
  $$
    \map{u_{g,h}}{f(gh)}{\rho(h\inv)(f(g))}
  $$
  satisfying the following two conditions: First, it is required that
$$u_{g,1}: f(g)\to \rho(1)(f(g))$$
coincides with $\phi_{1, f(g)}^{-1}$, see Definition \ref{2-rep-Def}(d). 
Second, it is required that for every $g\in G$ and every
  $h_1,h_2\in H$, the diagram

  \xymatrix{
  {f(gh_1h_2)} \ar[0,2]^{u_{gh_1,h_2}}\ar[2,0]_{u_{g,h_1h_2}}&&
  {\rho(h_1\inv)(f(gh_2))} \ar[2,0]^{\rho(h_2\inv)u_{g,h_1}}\\
  \\
  {\rho((h_1h_2)\inv)(f(g))}&& \rho(h_2\inv)\rho(h_1\inv)(f(g))
  \ar[0,-2]^{\phi_{h_1\inv,h_2\inv}}
  }
 commutes.

A morphism in $\ind HG(\mV)$ between two systems $(f, u=(u_{g,h}))$
and $(f', u'=(u'_{g,h}))$ is a system of morphisms $f(g)\to f'(g)$ in
$\mV$, given got each $g\in G$ and commuting with the $u_{g,h}$ and
$u'_{g,h}$. 

We define a left action $\sigma$ of  $G$  on $\ind HG(\mV)$ by
$$(\sigma(g_1)f)(g) = f(g_1^{-1}g), \quad (\sigma(g_1)u)_{g,h} = u_{g_1^{-1} g, h}.$$
\end{Def}

\begin{Rem} Consider the category $\prod_{g\in G}\mV$ whose objects are all
maps $G\to \ob \mV$. This category has a left $H$-action $\xi$ given by
$$(\xi(h) f)(g) = \rho(h)( f(gh)).$$
One sees immediately that 
$$\ind HG(\mV) = \bigl(\prod_{g\in G}\mV\bigr)^{H}$$
is identified with the category of $H$-equivariant objects in $\prod_{g\in G} \mV$. 
\end{Rem}

An explicit description of $\ind HG\mV$ is given as
follows (compare this to the classical definition of
the induced representation as in, e.g.,  \cite{Serre:representations}):

Let $m$ be the index of $H$ in $G$.
The underlying category of $\ind HG\mV$ is then identified with
 $\mV^m$. Such an identification is obtained by
  picking a
  system of representatives $$\mR = \{r_1,\dots r_m\}$$ of left cosets of
  $H$ in $G$ and associating to every map $f$ as above the
system $( f(r_1), ..., f(r_m))$.

 We view $\ind
  HG\rho(g)$ as $m\times m$ matrix whose
  entries are functors from $\mV$ to $\mV$. Then
  $$\(\ind HG\rho(g)\)_{ij}=
  \begin{cases}
    \rho(h) & \text{if}\quad gr_j = r_ih, h\in H\\
    0       & \text{else.}
  \end{cases}
  $$
  Note that in each row and each column, there is exactly one block
  entry, and that therefore such a matrix gives a functor from $\mV^m$
  to $\mV^m$.

\noindent  {\bf Composition} is defined as follows:
$$
\(\ind HG\rho(g_1)\)\circ_1\(\ind HG\rho(g_2)\)_{ik} =
\phantom{XXXXXXXXXXXXXXXX}$$
$$ {\phantom{XXXXXX} = \begin{cases}
    \rho(h_1)\circ_1\rho(h_2) & \text{if}\quad g_1r_j = r_ih_1 \text{
      and }  g_2r_k = r_jh_2 \\
    0       & \text{else.}
  \end{cases}}
$$
  In the case that this is not zero,
  $$\(\ind HG\rho(g_1g_2)\)_{ik} = \rho(h_1h_2),$$
  since in this case
  $$g_1g_2r_k = r_ih_1h_2.$$
  On this block the composition isomorphism is given by the
  2-isomorphism
  $$\rho(h_1)\circ_1\rho(h_2)\Rightarrow\rho(h_1h_2).$$
  Similarly, the isomorphism
  $$\ind HG\rho(1)\Rightarrow {1_{\mV^m}}$$
  is given by the corresponding map for $\rho$.
\begin{Prop}[{compare \cite[Ex.3.4.,Th.2]{Ostrik} and \cite[6.5]{Elgueta:2groups}}]
  Let $G$ be a finite group, and let $\mV=\Vect^{\oplus n}_k$. Then
  any \trep $\rho$ of $G$ in $\mV$ is isomorphic to a direct sum
  $$
    \bigoplus_{i=1}^m\ind{H_i}G\rho_{\omega_i},
  $$
  where the $H_i$ are subgroups of $G$, $\omega_i\in H^2(H_i,k^*)$, and
  $\rho_{\omega_i}$ is the $1$-dimensional \trep corresponding to
  $\omega_i$. Moreover, the system of $(H_i,\omega_i)$ is determined by
  the $G$-action on $\mV$ uniquely up to conjugation.
\end{Prop}
\begin{Pf}{}
  By Lemma \ref{perm-Lem}, $\rho$ defines a homomorphism from $G$ to $\Sn$,
  i.e., a $G$-action on the set $\{1,\dots,n\}$. Let $O_1,\dots,O_m$ be the
  orbits of this action. It follows that, after renumbering of $1,\dots, n$,
  the $2$-matrices $\rho(g)$ become block diagonal with blocks of sizes
  $|O_1|,\dots,|O_m|$. Hence
  $$
    \rho\cong\bigoplus_{i=1}^m \rho_i.
  $$
  Let $H_i$ be the stabilizer of an element of $O_i$. For $h\in H_i$, the
  $2$-matrix $\rho_i(h)$ is diagonal with the same $1$-dimensional vector space
  $L_i(h)$ on the diagonal. We conclude that
  $$
    \rho_i\cong \ind{H_i}G(\rho_{\omega_i}),
  $$
  where $\rho_{\omega_i}$ is the $1$-dimensional \trep of $H_i$ corresponding
  to the system
  $(L_i(h),h\in H_i)$.
\end{Pf}
\begin{Rem}\label{Ost}
  It follows that in the particular case of the proposition, the approach of
  Ostrik is equivalent to ours. In fact, the category
  $$
    \(\ind HG\(\rho_\omega\)\)^G, \quad H\sub G, \quad \omega\in H^2(H,k^*),
  $$
  is identified, as a $\mathcal {R}ep(G)$-module category, with the category of
  projective representations of $H$ with the central charge $\omega$.
\end{Rem}

\subsection{The character of the induced $2$-representation}
The aim of this subsection is to prove the following theorem.
\begin{Thm}\label{ind-char-Thm}
  Let $\mV$ be $k$-linear.
  The categorical trace $\mathbb Tr$ takes induced \treps into induced
  representations of groupoids. That is,
  $$
    \ttr{\ind HG\rho}\cong\ind{\Lambda(H)}\LG(\ttr\rho),
  $$
  as representations of $\LG$.
\end{Thm}
\begin{Cor}[{compare \cite[Thm D]{Hopkins:Kuhn:Ravenel}}]
  Assume that the $\ttr{\rho(h)}$ are finite dimensional.
  Let $\chi$ denote the 2-character of $\rho$.
  The 2-character of the induced representation is given by
  $$
    \chi_{\on{ind}}(g,h) = \GG H \sum_{s^{-1}(g,h)s\in H\times
      H}\chi(s^{-1}gs,s^{-1}hs).
  $$
\end{Cor}
\begin{Pf}{{of Theorem \ref{ind-char-Thm}}}
  We want to compute
  $$
    \Xind := \ttr{\ind HG\rho}
  $$
  as representation of
  $$
    \LG\simeq\coprod_{[g]_G}C_G(g).
  $$
  Let $\mR$ be a system of representatives of $G/H$, and fix
  $g\in G$. The underlying vector space of $\Xind(g)$ is the sum over
  all $r\in\mR$ which produce a diagonal block entry in
  $\on{ind}_\rho(g)$,
  \begin{equation}
    \label{Xind-Eqn}
    \Xind(g) = \bigoplus_{r^{-1}gr\in H}\ttr{\rho(r^{-1}gr)},
  \end{equation}
  (compare \cite{Serre:representations}).
  We need to determine the action of $C_G(g)$ on $\Xind(g)$.
  For this purpose, we replace our system of representatives $\mR$ in
  a convenient way:
  The decomposition
  $$
    [g]_G\cap H = [h_1]_H\cup\dots\cup [h_l]_H
  $$
  induces a decomposition
  $$
    \{r\in\mR\mid \con rg\in H \}= \bigcup_{i=1}^{l}\mR_i,
  $$
  with
  $$
    \mR_i := \{r\in\mR\mid\con rg\in[h_i]_H\}.
  $$
  We fix $i$, pick $r_i\in\mR_i$, and write $h_i:=\con {r_i}g$.
  \begin{Lem}
    \label{representatives-Lem}
    We can replace the elements of $\mR_i$ in such a way that left
    multiplication with $r_i^{-1}$ maps $\mR_i$ bijectively into a system of
    representatives of $$C_G(h_i)/C_H(h_i).$$
  \end{Lem}
  \begin{Pf}{}
    If $r\in\mR_i$ satisfies
    $$
      \con rg=\con h{h_i},
    $$
    we replace $r$ by $rh^{-1}$, which represents the same left coset
    of $G/H$ as $r$ does.
    Note that
    \begin{equation}
      \label{representatives-Eqn}
      (rh\inv)\inv grh\inv = h_i.
    \end{equation}
    We have
    $$\con{(r_i^{-1}rh^{-1})}{h_i} = h_i,$$
    therefore $r_i\inv (rh\inv)\in C_G(h_i)$.
    Assume now that we have replaced $\mR_i$ in this way.
    To prove that left multiplication with $r_i\inv$ is injective, let
    $r\neq r' \in\mR_i$. Then
    $$(r_i^{-1}r)^{-1}r_i^{-1}r' = r^{-1}r'$$ is not in $H$, and
    therefore $r_i^{-1}r'$ and $r_i^{-1}r$ are in different left
    cosets of $C_H(h_i)$ in $C_G(h_i)$. To prove
    surjectivity, let $\tilde g\in C_G(h_i)$. Write
    $$r_i\tilde g =rh$$
    with $r\in\mR$ and $h\in H$. Then
    $$
      \con rg = h\con{\tilde g}{\con {r_i}g}h^{-1}=h\con{\tilde
        g}{h_i}h^{-1}= hh_ih^{-1}.
    $$
    Therefore,
    $r$ is in $\mR_i$, and it follows from the identity
    (\ref{representatives-Eqn}) that $$\con rg=h_i.$$
    Thus $r_i^{-1}r=\tilde gh$ is in the same left coset of $C_H(h_i)$
    in $C_G(h_i)$ as $\tilde g$ is.
  \end{Pf}
  Let $\alpha_i$ denote the composition
  $$C_H(h_i)\hookrightarrow C_G(h_i){\to}C_G(g),$$
  where the second map is conjugation by $r_i^{-1}$.
  Recall that as representation of $C_G(g)$,
  \begin{eqnarray}\label{ttr-ind-Eqn}
    \(\ind{\Lambda(H)}\LG\pi\)(g) & = & \bigoplus_{i=1}^{l}
           \ind{\alpha_i}{}\pi(h_i).
  \end{eqnarray}
  \begin{Lem}\label{aaa}
    As a representation of $C_G$,
    $$
      \Xind(g)\cong\bigoplus_{i=1}^{l}\ind{\alpha_i}{}\ttr{\rho(h_i)}.
    $$
  \end{Lem}
  \begin{Pf}{} Let $f\in C_G(g)$, and let $r\in\mR_i$. Write
    $$
      fr=\tilde rh,
    $$
    with $\tilde r\in\mR$ and $h\in H$.
    We claim that $\tilde r$ is also in $\mR_i$ and that $h$ is in
    $C_H(h_i)$. This follows from
    $$
      \con{\tilde r}g =
      h\con{r}{\con fg}h^{-1} = h\con{r}{g}h^{-1} = hh_ih^{-1} = h_i
    $$
    as in the proof of Lemma \ref{representatives-Lem}.
  We are now ready to compute the block entry corresponding to $(r,r)$
  of
  $$
    \on{ind}_\rho(f^{-1})\circ_1 \on{ind}_\rho(g)\circ_1
    \on{ind}_\rho(f):
  $$
  \begin{eqnarray*}
    fr = \tilde rh &\text{gives}& (\on{ind}_\rho(f))_{\tilde rr} =
    \rho(h),\\
    g\tilde r = \tilde rh_i &\text{gives}& (\on{ind}_\rho(g))_{\tilde
      r\tilde r} = \rho(h_i),\\
    f^{-1}\tilde r = rh^{-1} &\text{gives}&
    (\on{ind}_\rho(f^{-1}))_{r\tilde r} = \rho(h^{-1})
  \end{eqnarray*}
  and all other block entries in these rows and columns are zero. Thus
  \begin{equation}
    \label{ind-conj-Eqn}
    \left(\on{ind}_\rho(f^{-1})\circ_1 \on{ind}_\rho(g)\circ_1
    \on{ind}_\rho(f)\right)_{rr} =
    \rho(h^{-1})\circ_1 \rho(h_i)\circ_1\rho(h),
  \end{equation}
  and the $2$-morphism from (\ref{ind-conj-Eqn}) to
  $$\left(\on{ind}_\rho(g)\right)_{rr} = \rho(h_i)$$
  is the $2$-morphism
  $$
    \rho(h^{-1})\circ_1 \rho(h_i)\circ_1\rho(h)\Rightarrow\rho(h_i)
  $$
  corresponding to $h^{-1}h_ih=h_i$.
  This proves that the action of $C_G(g)$ on $\Xind(g)$ decomposes into
  actions on
  $$
    \bigoplus_{r\in\mR_i}\rho(h_i).
  $$
  More precisely, if $fr=\tilde rh$, then $f$ maps the summand
  corresponding to $r$ to the one corresponding to $\tilde r$ by
  $$
    \map{h}{\rho(h_i)}{\rho(h_i)}.
  $$
  But
  $$
    fr = \tilde rh\iff(\con{r_i}{f})(r_i^{-1}r) = (r_i^{-1}\tilde r)h,
  $$
  and the action of $$\con{r_i}{f}\in C_G(h_i)$$ on
  $$\ind{C_H(h_i)}{C_G(h_i)}\rho(h_i)$$ is given by
  $$
    \map{h}{r_i^{-1}r\rho(h_i)}{r_i^{-1}\tilde r\rho(h_i)}.
  $$
Lemma \ref{aaa} is proved.
\end{Pf}
This now completes the proof of Theorem \ref{ind-char-Thm}.
\end{Pf}
\section{Some further questions}
\label{comp-HKR}
\subsection{Inertia orbifolds} Recall that a {\em Lie groupoid}
is a groupoid $\Gamma $ enriched in the category of $C^\infty$-manifolds,
i.e., such that $\ob (\Gamma)$ and $\mor (\Gamma)$ are $C^\infty$-manifolds and
all the structure maps (composition, inverses, units) are smooth. 
See  \cite{Mackenzie} for more details. An {\em orbifold}, cf.
\cite{Moerdijk}, is a Lie groupoid $G$ such that all stabilizer groups
$$Hom_{\Gamma} (x,x), \quad x\in\ob (\Gamma),$$
are finite. The construction of an inertia groupoid $\Lambda (\Gamma)$, see
Definition \ref{in-group}, can be applied to a Lie groupoid (resp. orbifold) $\Gamma$
and the result is again a Lie groupoid (resp. orbifold). 

\begin{Exa}[Global quotient groupoids] Let $M$ be a manifold and $G$ be a Lie
group acting on $M$. Then we have a Lie groupoid  $M\mmod G$ with
$$\ob(M\mmod G) = M, \quad Hom_{M\mmod G}(x,y) = \{ g\in G: g(x) = y\}.$$
Thus $\mor (M\mmod G) = G\times M$. 
If the stabilizer of each $x\in M$ is finite, then $M\mmod G$ is
an orbifold, known as the {\em global quotient orbifold}. 
The latter condition is automatically satisfied, if $G$ itself is finite. 
In this case the  inertia orbifold of $M\mmod G$  can be identified as follows:
$$
  \Lambda(M\mmod G) = \coprod_{[g]_G} M^g\mmod C_G(g).
$$
Here the disjoint union is over the conjugacy classes in $G$, 
and $M^g$ stands for the $g$-fixed point locus of $M$.
\end{Exa}

Recall further that equivariant $K$-theory of a manifold with a finite group
action is a particular case of a more general concept of orbifold $K$-theory
$K_\orb(\Gamma)$ 
defined for any orbifold $\Gamma$. This particular case 
corresponds to  a global quotient orbifold:
$$
  K_\orb(M\mmod G) \cong K_G(M).
$$
Our $2$-character map
\begin{equation}\label{trace-Eq}
  \map{\mathbb{T}r}{2\on{Rep}(G)}{\on{Rep}(\Lambda(G))}
\end{equation}
should be compared to the orbifold Chern character map defined by
Adem-Ruan and interpreted by Moerdijk \cite[p. 18]{Moerdijk} as a map
\begin{equation}\label{Chern-Eq}
K^\bullet(\Gamma)\otimes \mathbb{C}\to \prod_i H^{2i+\bullet}(\Lambda(\Gamma),
\mathbb{C}).
\end{equation}
Here  $\Gamma$ is any orbifold whose  quotient space 
(i.e., the space of isomorphism classes of objects) is compact.

This suggests that
\eqref{trace-Eq} has a generalization for an arbitrary orbifold $\Gamma$ as
above, yielding a transformation
$$2K_{\on{orb}}(\Gamma)\to K_{\on{orb}}(\Lambda(\Gamma)).$$
Here $2K_{\on{\orb}}(\Gamma)$ is a (yet to be defined) orbifold/equivariant
version of the 2-vector bundle $K$-theory of \cite{Baas:Dundas:Rognes}.
Recall that the non-orbifold $2K$ is interpreted as some approximation
to elliptic cohomology. Therefore the orbifold version is to be regarded
as a geometric version of equivariant elliptic cohomology, thus making
more precise our point in the introduction.
Note that inertia orbifolds also turn up in the original paper \cite{Hopkins:Kuhn:Ravenel},
where working at chromatic level $n$ requires using $n$-fold iterated inertia
orbifolds $\Lambda^n(\Gamma)$.

\vskip .2cm

\subsection{The Todd genus of $X^g$}
Let $G$ be a finite group 
acting on a compact $d$-dimensional complex manifold $X$. For
  $g\in G$, the fixed point locus $X^g$ is then a compact complex
  submanifold. Consider the graded vector spaces
  $$
    \mathbb T(g):=H^\bullet(X^g,\mathcal{ O}),
  $$
where $\mathcal{O}$ is the sheaf of holomorphic functions on $X^g$.
  Clearly, these are conjugation equivariant, i.e., they form a
  representation of $\LG$. Its character is a $2$-class function
  $$\chi_X\in 2\operatorname{Cl}(G; \ZZ[\zeta_N]), \quad
    \chi_X(g,h) = \tr\bigl\{ h_*: \mathbb T(g)\to \mathbb{T}(g)\bigr\}.
  $$
  This function takes values in the cyclotomic ring $\ZZ[\zeta_N]$,
  where $\zeta_N$ is an $N^{th}$ root of $1$ and $N$ is the order of
  $G$.

On the other hand, the Hopkins-Kuhn-Ravenel theory
\cite{Hopkins:Kuhn:Ravenel} also provides a 2-class function
associated with the $G$-manifold $X$. 
  Let $[X]\in\MU(\BG)$ denote the image of the equivariant cobordism
  class of $X$. Fix a prime $p$, and let $E=E_2$ be the second Morava
  $E$-theory at $p$. Recall that $E$ comes with a canonical natural
  transformation of cohomology theories
  $$\map\phi{\MU^*(-)}{E^*(-)}.$$
Let now $G$ be finite. In this situation Hopkins, Kuhn and Ravenel
constructed a map 
$$\alpha: E^*(BG)\to 2\operatorname{Cl}(G;D),$$
for a certain ring
$$
  D=\Colim_{n}D_n,
$$
where $D_n$ is  known as the ring of Drinfeld level $p^n$
structures on the formal group $E^*(\pt)$:
$$
  D_n := E^0(B(\ZZ/p^n\ZZ)^2)/\text{(annihilators of nontrivial Euler classes)}.
$$
The ring $D_n$ can be seen as the second chromatic analog of the cyclotomic ring
$\ZZ[\zeta_{p^n}]$ which corresponds to level $p^n$ structures on the multiplicative
group. In fact, a version of the Weil pairing \cite{AS} shows that $D_n$ contains
$\ZZ[\zeta_{p^n}]$. The 2-class function $a(y)$, $y\in E^*(BG)$ is defined by
$$a(y)(g,h):=(g,h)_n^*(y)\in D_n$$ 
where $(g,h)$ is a pair of commuting $p$-power order elements of $G$ and
$$(g,h)_n: (\ZZ/p^n\ZZ)^2\to G, \quad p^n= \max( \operatorname{ord}(g),
\operatorname{ord}(h)),$$
is the  homomorphism corresponding to $(g,h)$.

{\bf Question 1:} Is there a natural $2$-representation in some category of
sheaves associated to $X$ whose categorical character is  $\mathbb T$?

{\bf Question 2:} What is the relationship between $a(\phi[X])(g,h)$ 
and $\chi_X(g,h)$?


\end{document}